\documentclass[12pt,reqno]{amsart}

\usepackage{amssymb,amsmath,amsbsy,eucal,color,graphicx,subfigure,mathrsfs,latexsym}

\numberwithin{equation}{section}

\newtheorem{theorem}{Theorem}[section]

\theoremstyle{definition}

\theoremstyle{remark}
\newtheorem{remark}[theorem]{Remark}
\newtheorem{assumption}{Assumption}[section]

\numberwithin{equation}{section}

\newcommand{\vn}{{\vec{n}}}

\newcommand{\Ct}{{\mathcal T}}

\begin{document}

\title[ The stability in $L^{\infty}$ and $W^{1,\infty}$ in nonconvex polyhedra]{
The pointwise stabilities of piecewise linear finite element method on non-obtuse tetrahedral meshes of nonconvex polyhedra}

\author{Huadong Gao}
\address{School of Mathematics and Statistics, 
Huazhong University of Science and Technology, 
Wuhan 430074, China} 
\email{huadong@hust.edu.cn}

\author{Weifeng Qiu}
\address{Department of Mathematics,
City University of Hong Kong, Kowloon, Hong Kong, China}
\email{weifeqiu@cityu.edu.hk}

\thanks{The second author is the corresponding author.}

\begin{abstract} 
Let $\Omega$ be a Lipschitz polyhedral (can be nonconvex) domain in $\mathbb{R}^{3}$, and 
$V_{h}$ denotes the finite element space  of continuous piecewise linear polynomials. 
On non-obtuse quasi-uniform tetrahedral meshes, we prove that 
the finite element projection $R_{h}u$  of $u \in H^{1}(\Omega) \cap C(\overline{\Omega})$ 
(with $R_{h} u$ interpolating $u$ at the boundary nodes) satisfies 
\begin{align*}
\Vert R_{h} u\Vert_{L^{\infty}(\Omega)} \leq C \vert \log h \vert \Vert u\Vert_{L^{\infty}(\Omega)}.
\end{align*}
If we further assume $u \in W^{1,\infty}(\Omega)$, then 
\begin{align*}
\Vert R_{h} u\Vert_{W^{1, \infty}(\Omega)} \leq C \vert \log h \vert \Vert u\Vert_{W^{1, \infty}(\Omega)}.
\end{align*}
\end{abstract}



\maketitle

\section{Introduction}

In this paper we consider the the Ritz projection $R_{h} u \in V_{r, h}$ of $u \in H^{1}(\Omega) 
\cap C(\overline{\Omega})$  satisfying 
\begin{align}
\label{def_R_h_general}
(\nabla R_{h}u, \nabla v_{h})_{\Omega} = (\nabla u, \nabla v_{h})_{\Omega}, \quad \forall v_{h} \in V_{r, h}^{0},
\end{align} 
where $V_{r, h}$ is the finite element subspace of $H^{1}(\Omega)$ composed of piecewise polynomials 
of degree $r$ ($r \geq 1$), $V_{r, h}^{0} = H^{1}_{0}(\Omega) \cap V_{r, h}$, 
and $R_{h}u$ interpolates $u$ at the boundary nodes on $\partial\Omega$. 
In fact, $R_{h}u$ is the finite element projection of $u$ onto $V_{r, h}$ for the model problem 
\begin{align}
\label{poisson_eq}
\Delta u = f \text{ in } \Omega,
\end{align} 
with Dirichlet boundary condition on $\partial\Omega$. 

Our motivation is to establish the stability in $L^{\infty}(\Omega)$
\begin{subequations}
\label{stabs_L_infty}
\begin{align}
\label{stab1_L_infty}
& \Vert R_{h} u \Vert_{L^{\infty}(\Omega)} \leq C \Vert u\Vert_{L^{\infty}(\Omega)}, \\
\label{stab2_L_infty}
& \text{or } \Vert R_{h} u \Vert_{L^{\infty}(\Omega)} \leq C \vert \log h \vert  \Vert u\Vert_{L^{\infty}(\Omega)}; 
\end{align}
\end{subequations}
and the stability in $W^{1,\infty}(\Omega)$ (if $u \in W^{1,\infty}(\Omega)$) 
\begin{subequations}
\label{stabs_W_1_infty}
\begin{align}
\label{stab1_W_1_infty}
& \Vert R_{h} u \Vert_{W^{1,\infty}(\Omega)} \leq C \Vert u\Vert_{W^{1,\infty}(\Omega)}, \\ 
\label{stab2_W_1_infty}
& \text{or } \Vert R_{h} u \Vert_{W^{1,\infty}(\Omega)} \leq C \vert \log h \vert  \Vert u\Vert_{W^{1,\infty}(\Omega)}.
\end{align}
\end{subequations}

There are a lot of important works for estimates (\ref{stabs_L_infty}) and (\ref{stabs_W_1_infty}).   
\cite{Natterer1975} and \cite{Scott1976} are the first contributions for general quasi-uniform meshes. 
On convex polygonal domains,  \cite{Natterer1975} considered piecewise linear ($r = 1$) approximation 
while \cite{Scott1976} treated the finite element approximation (for any $r\geq 1$) to Neumann problem of 
(\ref{poisson_eq}).  \cite{Schatz1980} proved (\ref{stabs_L_infty}) and (\ref{stabs_W_1_infty}) on polygonal 
(can be nonconvex)  domains.  When $r=1$,  the estimates provided in \cite{Schatz1980} are (\ref{stab2_L_infty}) 
and (\ref{stab2_W_1_infty}).  Thus,  estimates (\ref{stabs_L_infty}) and (\ref{stabs_W_1_infty}) are valid for 
most practical domains in $\mathbb{R}^{2}$. On the contrast,  in three dimensional space,  all existing works 
\cite{BS,Demlow2012, Guzman09, LeykekhmanLi2021, LeykekhmanVexler2016, Rannacher1976, RannacherScott1982,  
Schatz1998, SchatzWahlbinq982} for estimates (\ref{stabs_L_infty}) and (\ref{stabs_W_1_infty}) are available 
on either domains with smooth boundary or convex polyhedral domains (Instead of explicit assumptions on domains, 
\cite{BS} needs $\Vert w \Vert_{W^{2,p}(\Omega)} \leq C \Vert \Delta w \Vert_{L^{p}(\Omega)}$ 
for some $p > 3$ in three dimensional space,  for any function $w$ with zero trace on $\partial\Omega$).     

In this paper, we prove that if the meshes are non-obtuse (all internal dihedral angles of all tetrahedral elements 
are less than or equal to $\frac{\pi}{2}$), then estimates (\ref{stab2_L_infty}) and (\ref{stab2_W_1_infty}) 
hold for the finite element projection (\ref{def_R_h_general}) with piecewise linear finite element space 
$V_{h} = V_{1,h}$ ($r=1$). 

In Section $2$, we provide the main results and all assumptions. In Section $3$, we show the proofs of our main results.

\section{Main results}

Let $\Omega$ be a Lipschitz polyhedra (can be nonconvex) in $\mathbb{R}^{3}$. 
We denote by $\Ct_{h}$ quasi-uniform conforming tetrahedral meshes of $\Omega$. 
We define $V_{h} = H^{1}(\Omega) \cap P_{1}(\Ct_{h})$ and  $V_{h}^{0} = H_{0}^{1}(\Omega) \cap V_{h}$.  
For any $u \in H^{1}(\Omega)$,  we introduce the Ritz projection $R_{h}u \in V_{h}$ to satisfy 
\begin{align}
\label{def_R_h}
(\nabla R_{h}u, \nabla v_{h})_{\Omega} = (\nabla u, \nabla v_{h})_{\Omega}, \quad \forall v_{h} \in V_{h}^{0},
\end{align}
where $R_{h}u$ interpolates $u$ at the boundary nodes on $\partial\Omega$.  In fact, $R_{h}u$ is the finite element 
projection of $u$ onto $V_{h}$, and (\ref{def_R_h}) is exactly the finite element projection (\ref{def_R_h_general}) 
with $r = 1$.
 
\begin{assumption}
\label{assump1}
For any $T \in \Ct_{h}$, all internal dihedral angles of the tetrahedral element $T$ are less than or equal to $\frac{\pi}{2}$.  
$\Ct_{h}$ is called non-obtuse tetrahedral meshes of $\Omega$.
\end{assumption}

\begin{assumption}
\label{assump2}
The mesh $\Ct_{h}$ of $\Omega$ can be extended to a larger convex domain $\tilde{\Omega}$ quasi-uniformly with 
$\Omega \Subset  \tilde{\Omega}$. We denote by $\tilde{\Ct}_{h}$ the extension of $\Ct_{h}$ on $\tilde{\Omega}$.  
\end{assumption}

\begin{remark}
We don't require $\tilde{\Ct}_{h}$ introduced in Assumption~\ref{assump2} to be non-obtuse for 
all tetrahedral elements.  Only elements $T \in \Ct_{h}$ need to be non-obtuse. 
\end{remark}

\begin{theorem}
\label{thm_L_infty}
If Assumption (\ref{assump1}) and Assumption (\ref{assump2}) hold,  then there is a positive constant $C$ such that 
for any $u \in H^{1}(\Omega) \cap C(\overline{\Omega})$, 
\begin{align*}
\Vert R_{h}u \Vert_{L^{\infty}(\Omega)} \leq C \vert \log h \vert \Vert u\Vert_{L^{\infty}(\Omega)}.
\end{align*}
\end{theorem}

\begin{theorem}
\label{thm_W_1_infty}
If Assumption (\ref{assump1}) and Assumption (\ref{assump2}) hold, then there is a positive constant $C$ such that 
for any $u \in W^{1,\infty}(\Omega)$, 
\begin{align*}
\Vert R_{h}u \Vert_{W^{1, \infty}(\Omega)} \leq C \vert \log h \vert \Vert u\Vert_{W^{1, \infty}(\Omega)}.
\end{align*}
\end{theorem}

\section{Analysis}

\begin{proof}
(Proof of Theorem~\ref{thm_L_infty}) 
Since $u \in C(\overline{\Omega})$,  we denote by $\tilde{u}$ the extension of $u$ to $\tilde{\Omega}$,  such that 
$u\in C_{0}(\tilde{\Omega})$ and $\Vert \tilde{u} \Vert_{L^{\infty}(\tilde{\Omega})} = \Vert u \Vert_{L^{\infty}(\Omega)}$. 
The existence of $\tilde{u}$ satisfying the above two properties follows from the facts that $u \in C(\overline{\Omega})$ and 
the Whitney type extension operator $\mathcal{E}_{0}$ in Section $2.2$ of Chapter $6$ in \cite{Stein70} (see ($8$) 
and the proposition in Section $2.2$ of Chapter $6$ in \cite{Stein70}). We would like to emphasize that we don't need 
$\tilde{u} \in H^{1}(\tilde{\Omega})$.


We define $\tilde{V}_{h}^{0} = H_{0}^{1}(\tilde{\Omega}) \cap P_{1}(\tilde{\Ct}_{h})$. 
Let $\tilde{u}_{h} \in \tilde{V}_{h}^{0}$ satisfy 
\begin{align}
\label{def_tilde_u_h}
(\nabla \tilde{u}_{h}, \nabla \tilde{v}_{h})_{\tilde{\Omega}} = 
\Sigma_{T \in \tilde{\Ct}_{h}} \left(  - (\tilde{u}, \Delta \tilde{v}_{h})_{T} + \langle \tilde{u}, 
\nabla \tilde{v}_{h} \cdot \vn \rangle_{\partial T} \right), 
\quad \forall \tilde{v}_{h} \in \tilde{V}_{h}^{0}.
\end{align}
Here $\vn$ is the outward unit normal vector along $\partial T$ for any $T\in \Ct_{h}$. 
For any $v_{h} \in V_{h}^{0}=H_{0}^{1}(\Omega)\cap P_{1}(\Ct_{h})$, we denote by $\tilde{v}_{h} \in \tilde{V}_{h}^{0}$ 
the zero extension of $v_{h}$ to $\tilde{\Omega}$. 
By (\ref{def_tilde_u_h}) and the definition of $\tilde{v}_{h}$,  it is easy to see that 
\begin{align}
\label{tilde_u_h}
& (\nabla \tilde{u}_{h}, \nabla v_{h})_{\Omega}
 =  (\nabla \tilde{u}_{h}, \nabla \tilde{v}_{h})_{\tilde{\Omega}} \\ 
\nonumber 
= & \Sigma_{T \in \tilde{\Ct}_{h}} \left(  - (\tilde{u}, \Delta \tilde{v}_{h})_{T} + \langle \tilde{u}, 
\nabla \tilde{v}_{h} \cdot \vn \rangle_{\partial T} \right)  \\ 
\nonumber
= & \Sigma_{T \in \Ct_{h}} \left(  - (\tilde{u}, \Delta \tilde{v}_{h})_{T} + \langle \tilde{u}, 
\nabla \tilde{v}_{h} \cdot \vn \rangle_{\partial T} \right)  \\
\nonumber
= & \Sigma_{T \in \Ct_{h}} \left(  - (u, \Delta v_{h})_{T} + \langle u, 
\nabla v_{h} \cdot \vn \rangle_{\partial T} \right) 
= (\nabla u, \nabla v_{h})_{\Omega}. 
\end{align}
The last equality holds since $u \in H^{1}(\Omega)$. 
On the other hand,  since $\tilde{\Omega}$ is convex and $\tilde{u} \in C_{0}(\tilde{\Omega})$, 
(\ref{def_tilde_u_h}) and \cite[Theorem~$12$]{LeykekhmanVexler2016} imply that 
\begin{align}
\label{convex_L_infty}
\Vert \tilde{u}_{h} \Vert_{L^{\infty}(\tilde{\Omega})} 
\leq C \vert \log h \vert \Vert \tilde{u} \Vert_{L^{\infty}(\tilde{\Omega})} 
= C \vert \log h \vert \Vert u\Vert_{L^{\infty}(\Omega)} .
\end{align}

We notice that $R_{h}u \in V_{h}=H^{1}(\Omega) \cap P_{1}(\Ct_{h})$ satisfies
\begin{align*}
(\nabla R_{h}u, \nabla v_{h})_{\Omega} = (\nabla u, \nabla v_{h})_{\Omega}, \quad \forall v_{h} \in V_{h}^{0}
= H_{0}^{1}(\Omega) \cap V_{h}.
\end{align*}
Thus, by the above equation and (\ref{tilde_u_h}), we have that $\left( R_{h}u - \tilde{u}_{h} \right) |_{\Omega} 
\in V_{h}$ and 
\begin{align*}
(\nabla (R_{h}u - \tilde{u}_{h}), \nabla v_{h})_{\Omega} = 0, \quad \forall v_{h} \in V_{h}^{0} 
= H_{0}^{1}(\Omega) \cap V_{h}.
\end{align*}
By Assumption~(\ref{assump1}) and \cite[Theorem~$3.2$ and Lemma~$5.1(\text{iii})$]{WangZhang2012} 
(or by \cite{CiarletRaviart1973, KorotovKrizekNeittaanmaki2000, KrizekQun1995}), the above equation implies that  
\begin{align}
\label{weak_max_principle}
\Vert R_{h}u - \tilde{u}_{h} \Vert_{L^{\infty}(\Omega)} \leq \Vert R_{h}u - \tilde{u}_{h}\Vert_{L^{\infty}(\partial \Omega)} 
\leq \Vert u\Vert_{L^{\infty}(\partial\Omega)} +  \Vert \tilde{u}_{h} \Vert_{L^{\infty}(\partial \Omega)}.
\end{align}

Thus, by (\ref{convex_L_infty}) and (\ref{weak_max_principle}),  it is easy to see that 
\begin{align*}
& \Vert R_{h} u \Vert_{L^{\infty}(\Omega)} 
\leq \Vert R_{h}u - \tilde{u}_{h} \Vert_{L^{\infty}(\Omega)}  + \Vert \tilde{u}_{h} \Vert_{L^{\infty}(\Omega)} \\
\leq & \Vert u\Vert_{L^{\infty}(\partial\Omega)} + \Vert \tilde{u}_{h} \Vert_{L^{\infty}(\partial \Omega)} +  \Vert \tilde{u}_{h} \Vert_{L^{\infty}(\Omega)} \\
\leq & \Vert u\Vert_{L^{\infty}(\Omega)} + 2  \Vert \tilde{u}_{h} \Vert_{L^{\infty}(\Omega)}
 \leq C \vert \log h \vert \Vert u\Vert_{L^{\infty}(\Omega)}.
\end{align*}

The proof is complete.
\end{proof}

\begin{proof}
(Proof of Theorem~\ref{thm_W_1_infty}) 
We denote by $I_{h}u$ the standard interpolation of $u$ on $V_{h} = H^{1}(\Omega) \cap P_{1}(\Ct_{h})$. 

By applying Theorem~\ref{thm_L_infty} to $u - I_{h}u$, we have 
\begin{align*}
\Vert R_{h}u - I_{h}u \Vert_{L^{\infty}(\Omega)} \leq C \vert \log h \vert \Vert u - I_{h}u \Vert_{L^{\infty}(\Omega)}.
\end{align*}

By inverse inequality and approximation properties of $I_{h}$, 
\begin{align*}
& \Vert R_{h} u \Vert_{W^{1,\infty}(\Omega)} \leq \Vert R_{h} u - I_{h} u \Vert_{W^{1,\infty}(\Omega)} 
+ \Vert I_{h}u\Vert_{W^{1,\infty}(\Omega)} \\
\leq & C h^{-1} \Vert R_{h}u - I_{h}u\Vert_{L^{\infty}(\Omega)} + C \Vert u\Vert_{W^{1,\infty}(\Omega)} 
\leq C \Vert u\Vert_{W^{1,\infty}(\Omega)}.
\end{align*}

The proof is complete.
\end{proof}

\section{Declarations}
\noindent {\bf Funding:}
Huadong Gao is partially supported by  National
Natural Science Foundation of China under grant number 11871234.
Weifeng Qiu is supported by 
a grant from the Research Grants Council of the Hong Kong Special Administrative Region, China 
(Project No. CityU 11302718).  \\

\noindent {\bf The Conflict of Interest Statement:}
No conflict of interest exists.\\

\noindent {\bf Availability of data and material:} 
Not applicable.\\

\noindent {\bf Code availability:}
Not applicable.\\

\noindent {\bf Authors' contributions:}
Huadong Gao and Weifeng Qiu have participated sufficiently in the work to take public
responsibility for the content, including participation in the concept, method,
analysis and writing. All authors certify that this material or similar material
has not been and will not be submitted to or published in any other publication.

\end{document}